\documentclass[12pt]{amsart}
\usepackage{amsmath,amstext, amsfonts,amssymb}
\usepackage{mathrsfs}
\usepackage[T2A]{fontenc}
\usepackage[utf8]{inputenc}
\usepackage[russian, english]{babel}
\usepackage{cite, color}
\usepackage{enumitem}
\usepackage{tikz-cd}

\newcommand{\RR}{\mathbb R}

\newtheorem{theorem}{Theorem}
\newtheorem{lemma}{Lemma}

\newtheorem{definition}{Definition}
\newtheorem{corollary}{Corollary}

\allowdisplaybreaks

\begin{document}

\title{Ostrowski-type inequalities in abstract distance spaces}

\author[V.~F.~Babenko]{Vladyslav Babenko}
\address{Department of Mathematical Analysis and Optimization, Oles Honchar Dnipro National University, Dnipro, Ukraine}

\author[V.~V.~Babenko]{Vira Babenko}
\address{Department of Mathematics and Computer Science, Drake University, Des Moines, USA}

\author[O.~V.~Kovalenko]{Oleg Kovalenko}
\address{Department of Mathematical Analysis and Optimization, Oles Honchar Dnipro National University, Dnipro, Ukraine}

\begin{abstract}
For non-empty sets X we define  notions of distance and pseudo metric with values in a partially ordered set that has a  smallest element $\theta $. If $h_X$ is a distance in $X$ (respectively, a pseudo metric in $X$), then the pair $(X,h_X)$ is called a distance (respectively, a pseudo metric) space. If $(T,h_T)$ and $(X,h_X)$ are pseudo metric spaces, $(Y,h_Y)$ is a distance space, and $H(T,X)$ is a class of Lipschitz mappings $f\colon T\to X$, for a broad family of mappings $\Lambda\colon H  (T,X)\to Y$, we obtain a sharp inequality that estimates the deviation $h_Y(\Lambda f(\cdot),\Lambda f(t))$ in terms of the function $h_T(\cdot, t)$. We also show that many known estimates of such kind are contained in our general result.
\end{abstract}
\maketitle

\section{Introduction}
Estimates for the deviation between the value of an operator $\Lambda$ at a function $f$ from some  class $\mathfrak{M}$ and the value of $\Lambda$ at some  depending on $f$ constant function from $\mathfrak{M}$ play an important role in approximation theory and numeric analysis. For example, estimates for the deviation of a value of a function $f\in \mathfrak{M}$ at some point from its mean value is of this kind. One of the first among such sharp estimates (where $\Lambda f=\frac 12 \int_{-1}^1 f(t)dt$) was obtained by Ostrowski~\cite{Ostrowski38}:
\begin{theorem} Let $f\colon [-1,1]\to \mathbb{R}$ be a differentiable function and let for all $t\in (-1,1)$, $|f'(t)|\leq 1$. Then for all $x\in [-1,1]$ the following inequality holds
\begin{equation}\label{widerclass}
\left|\frac 12 \int\limits_{-1}^1 f(t)dt - f(x)\right|
\leq \frac 12\int\limits_{-1}^1|t-x|dt=\frac{1+ x^2}{2}.
\end{equation}
The inequality is sharp in the sense that for each fixed $x\in [-1,1]$, the upper bound $\frac{1+ x^2}{2}$ cannot be reduced.
\end{theorem}

Note that inequality~\eqref{widerclass} holds for the wider class of functions that satisfy the Lipschitz condition
$$
|f(t_1)-f(t_2)|\le |t_1-t_2|\;\;\forall t_1,t_2\in [-1,1],    
$$
and becomes equality on the function $f_x(t)=|t-x|$, $t\in [-1,1]$.

Inequalities that estimate the deviation of a value of a function at some point from its mean value using some characteristics of the function, are sometimes called Ostrowski-type inequalities. Such inequalities were intensively studied, see for example~\cite{Dragomir17,Dragomir02}.

Let two non-empty sets   $X$ and $Y$ be given. The class $H(X,Y)$ of mappings $f\colon X\to Y$ that satisfy the Lipschitz condition can be defined in a standard way, if distances $h_X$ and  $h_Y$ are somehow defined in $X$ and $Y$:
\[
H(X,Y)=\{ f\colon X\to Y\colon h_Y(f(x_1),f(x_2))\le h_X(x_1,x_2) \;\;\forall x_1,x_2\in X\}.
\]

In this paper for a non-empty sets X we using a partially ordered set $M$ that has a smallest element, we define (see Section~\ref{s::notations}) a concept of a distance and a concept of a pseudo metric $h_X$ with values in $M$. Such distances (pseudo metrics) will be called abstract distances (abstract pseudo metrics) or $M$-distances ($M$-pseudo metrics). We define an $M$-distance space (an $M$-pseudo metric space) as a pair $(X,h_X)$ where $h_X$ is an $M$-distance on $X$ ($h_X$ is an $M$-pseudo metric 
on $X$ respectively). If $(T,h_T)$ and $(X,h_X)$ are $M$-pseudo metric spaces, $(Y,h_Y)$ is an $M$-distance space, for a broad family of mappings $\Lambda\colon H(T,X)\to Y$, we prove (see Section~\ref{s::classesAndInequality}) the main result of the paper -- a sharp inequality that estimates the deviation $h_Y(\Lambda f(\cdot),\Lambda f(t))$ (here $t\in T$ and $\Lambda f(t)$ is the value of $\Lambda$ at the function $\tau\to f(t)$, $\tau\in T$) in terms of some characteristic of the function $h_T(\cdot, t)$. We also show that many known estimates as well as many new estimates of such kind are contained in our general result, which is natural to be called an abstract version of the Ostrowski-type inequalities. 
In Section~\ref{s::Hw} we 
introduce a general concept of a modulus of continuity as a function $\omega \colon M\to M$ that satisfies certain conditions, as well as the corresponding classes $H^\omega (T,X)$. For a broad family of mappings $\Lambda\colon H^\omega (T,X)\to Y$, we obtain a sharp inequality that estimates the deviation $h_Y(\Lambda f(\cdot),\Lambda f(t))$ in terms of the value $\omega(h_T(\cdot,t))$.
Finally, in Section~\ref{s::agreementOfDistances}, we give rather general conditions for an  $M$-distance space $T$ to be an $M$-pseudo metric space.

\section{Notations and definitions}\label{s::notations}
The notion of a distance (in particular, a metric) plays an important role in many branches of mathematics. Definitions of numeric-valued distances or metrics and a detailed discussion of these notions can be found e.g., in monograph~\cite{Kirk}.  We refer to~\cite{HausdDistSpaces,fmDistSpaces,Jancovic,Zabrejko} for metrics that take value in more general sets. We consider a rather general definition for this notion.

A set $M$ with a reflexive, antisymmetric and transitive relation $\leq$ is called partially ordered.

\begin{definition}
Let $X$ be an arbitrary set and  $M$ be a partially ordered set that has a smallest element, which we denote by $\theta$ (i.e., $\theta\leq m$ for any $m\in M$). A function  $h_X\colon X\times X\to M$ is called an {\it $M$-distance} in $X$, if for arbitrary $x,y\in X$
\begin{enumerate}
\item $h_X(x,x)=\theta,$
\item $h_X(x,y)=h_X(y,x)$.
\end{enumerate}
The pair $(X,h_X)$ will be called an {\it $M$-distance space}.    
\end{definition}

In~\cite{HausdDistSpaces,fmDistSpaces} the notion of $M$-distance was introduced for the case when  $M$ is a partially ordered monoid.

Everywhere below, speaking of a partially ordered set $M$, we assume that some $M$-distance $h_M$ is defined in $M$.

\begin{definition}
We say that an $M$-distance $h_X$ in $X$ {\it agrees} with an $M$-distance $h_M$ in $M$, if 
\begin{equation}\label{metricsAgreementCond}
   h_M(h_X(x,x_1),h_X(x,x_2))
    \le 
    h_X(x_1,x_2)\;\;\;\forall x,x_1,x_2\in X. 
\end{equation}   
\end{definition}
 Note that inequality~\eqref{metricsAgreementCond} holds (and is equivalent to the triangle inequality) if $M=\RR_+$ with the usual metric, and $(X,h_X)$ is a pseudo metric space (for a definition of a pseudo metric and a pseudo metric space see, for example~\cite[Chapter~4]{Kelly}). That is why we introduce the following definition.

\begin{definition}
An $M$-distance $h_X$ on a set $X$ will be called an {\it $M$-pseudo metric}, if it agrees with $M$-distance $h_M$ i.e., inequality~\eqref{metricsAgreementCond} holds. In this case the pair  $(X,h_X)$ will be called an  {\it $M$-pseudo metric space}.    
\end{definition}

In Lemma~\ref{l::metricsByE} we will give a general sufficient condition
that an $M$-metric $h$ (see Section~\ref{s::agreementOfDistances} for the definition of $M$-metric) agrees with $h_M$.

We need the following lemma.
\begin{lemma}\label{l::Lemma1}
Let $T,X,Y$ be $M$-distance spaces. Then 
\begin{enumerate}
    \item If $f\in H(T,X)$ and $g\in H(X,Y)$, then $g\circ f\in H(T,Y)$. 
    \item If $f\in H(T,X)$ and $h_X$ is an $M$-pseudo metric, then one has $h_X(f(\cdot), f(t))\in H(T,M)$ for any fixed $t\in T$. In particular, if $T=X$, then $h_T(\cdot,t)\in H(T,M)$.  
 \end{enumerate}
\end{lemma}
\begin{proof}
The first statement of the lemma is obvious. 

If $f\in H(T,X)$ and $h_X$ is an $M$-pseudo metric, then for arbitrary $t_1,t_2\in T$,
$$
    h_M(h_X[f(t_1), f(t)], h_X[f(t_2),f(t)])
    \stackrel{\eqref{metricsAgreementCond}}{\le}
    h_X(f(t_1),f(t_2))
    \stackrel{f\in H(T,X)}{\le}
    h_T(t_1,t_2).
$$
Therefore $h_X(f(\cdot), f(t))\in H(T,M)$. If $T=X$ and $f(\tau)=\tau, \; \tau\in T,$ we obtain that $h_T(\cdot,t)\in H(T,M)$.
\end{proof}
\section{Classes of operators and an Ostrowski-type inequality}\label{s::classesAndInequality}
\begin{definition}
For  an $M$-distance space $X$, an operator  $\lambda\colon  H(X,M)\to M$ will be called {\it monotone}, if for arbitrary $u,v\in H(X,M)$ 
$$
        (\forall x\in X\; u(x)\le v(x))\implies (\lambda(u)\le \lambda(v)).
$$
\end{definition}
\begin{definition}
Let $T,Y$ be $M$-distance spaces,  X be an $M$-pseudo metric space, and $t\in T$ be  fixed. We say that an operator $\Lambda\colon H(T,X)\to Y$ and a monotone operator $\lambda\colon  H(T,M)\to M$ {\it agree}, if  $\forall f\in H(T,X)$ 
\begin{equation}\label{Agree}
h_Y(\Lambda f(\cdot),\Lambda f(t))\le 
\lambda( h_X(f(\cdot),f(t))).
\end{equation}
%{\color{red}In this case we write $\Lambda\ll\lambda$.}
\end{definition}
Here and everywhere below   $\Lambda f(t)$ means the value of the operator  $\Lambda$ on the constant function $\tau\mapsto f(t)$, $ \tau\in T$ (the same notation will be used for other operators whose arguments are functions). 

%{\color{blue} Може має смисл перенести сюди пояснення того що прикладами узгоджених операторів є інтеграл Бохнера і інтеграл Лебега, а також інтеграл від функції зі значеннями у $L$-просторі і інтеграл Лебега. Мені здпється що про Бохнера має сенс. А про решту можна щось сказати словами.}

\begin{theorem}\label{th::generalOstrowskiHw}
 Let $(T,h_T)$ and $(X,h_X)$ be $M$-pseudo metric spaces, $(Y,h_Y)$ be an $M$-distance space, and $t\in T$ be fixed. 
    Assume that an operator $\Lambda\colon H(T,X)\to Y$ and a monotone operator $\lambda\colon  H(T,M)\to M$ agree. Then for arbitrary function $f\in H(T,X)$ the following Ostrowski-type inequality holds:
    \begin{equation}\label{OstTypeInequality}
    h_Y(\Lambda f(\cdot), \Lambda f(t))\le 
\lambda (h_T(\cdot,t)).
    \end{equation}
    If
    \begin{equation}\label{lambdaProp}
        \lambda(\theta)=\theta,
    \end{equation} 
     and there exists  an operator $\phi_X\colon H(T,M)\to H(T,X)$ and $\phi_Y\in H(M,Y)$ with the following property
\begin{equation}\label{phiProp}
h_Y(\phi_Y(m),\phi_Y(\theta))=m\text{, if } m = \lambda(h_T(\cdot,t)),
    \end{equation}
such that the diagram
$$
\begin{tikzcd}
H(T,X) \arrow{r}{\Lambda}  & Y \\
H(T,M) \arrow{r}{\lambda} 
\arrow{u}{\phi_X}  & M  \arrow{u}[swap]{\phi_Y}
\end{tikzcd}    
$$
is commutative i.e., 
\begin{equation}\label{diagram}
\Lambda\circ \phi_X =\phi_Y\circ \lambda,
\end{equation}
then inequality~\eqref{OstTypeInequality} is sharp and becomes equality on the function
\begin{equation}\label{abstractOstrowskiHwExtremal}
f_t(\cdot)=\phi_{X}(h_T(\cdot,t)).
\end{equation}
\end{theorem}
\begin{proof}
 Let $f\in H(T,X)$. Since $(T,h_T)$ and $(X,h_X)$ are $M$-pseudo metric spaces, we have due to Lemma~\ref{l::Lemma1} that $h_X(f(\cdot),f(t))\in H(T,M)$ and $h_T(\cdot,t)\in H(T,M)$. So both of these functions belong to the domain of $\lambda$.

    Since the operators $\Lambda$ and $\lambda$ agree,  and the operator  $\lambda$ is monotone, for each  $f\in H(T,X)$ one has
    \[
    h_Y(\Lambda f(\cdot), \Lambda f(t))
    \le \lambda(h_X(f(\cdot),f(t)))
\le \lambda(h_T(\cdot,t)),
    \]
and inequality~\eqref{OstTypeInequality} is proved. 

The function from~\eqref{abstractOstrowskiHwExtremal} belongs to the class $H(T,X)$, since $h_T(\cdot,t)\in H(T,M)$ and $\phi_X\colon H(T,M)\to H(T,X)$. Using condition~\eqref{diagram}, one has
    \[
    h_Y(\Lambda f_t(\cdot), \Lambda f_t(t))=h_Y(\Lambda(\phi_X(h_T(\cdot,t))), \Lambda(\phi_X(h_T(t,t))))
    \]
    \[
    =h_Y((\Lambda\circ  \phi_X)(h_T(\cdot,t)), (\Lambda\circ\phi_X)(h_T(t,t)))
    \]
\[
\stackrel{\eqref{diagram}}{=}h_Y((\phi_Y\circ\lambda)(h_T(\cdot,t)),(\phi_Y\circ\lambda)(\theta))
\]
\[
=h_Y(\phi_Y(\lambda(h_T(\cdot,t))),\phi_Y(\lambda(\theta)))
\stackrel{\eqref{phiProp},\eqref{lambdaProp}}{=}
\lambda(h_T(\cdot,t)).
\]
The theorem is proved.
\end{proof}

Note that classes of operators that satisfy the properties analogues of  the properties~\eqref{Agree},~\eqref{phiProp} and~\eqref{diagram} were considered in~\cite{Kovalenko21a}.
In order to explain the nature of these properties we give the following example.
Condition~\eqref{Agree} is a relaxed version of the following condition:
$$h_Y(\Lambda f,\Lambda g)\le 
\lambda(h_X(f(\cdot),g(\cdot)))$$
for all $f,g$. If $M = \RR_+$, $X = Y$ is a Banach space, $\Lambda f = \int_{-1}^1 f(t)dt$ is the Bochner integral of $f$, and $\lambda$ is the Lebesgue integral on $[-1,1]$, then this condition becomes
$$
\left\|\int\limits_{-1}^1 f(t)dt - \int\limits_{-1}^1 g(t)dt\right\|\leq \int\limits_{-1}^1 \|f(t)-g(t)\|dt.
$$
Moreover, if $\alpha$ is an integrable real-valued function and $x\in X$, then 
$$
\int\limits_{-1}^1 \alpha(t)\cdot x dt = \left(\int\limits_{-1}^1\alpha(t)dt\right)\cdot x,
$$
i.e., condition~\eqref{diagram} is satisfied with $\phi_X$ and $\phi_Y$ being multiplication by a fixed element $x\in X$. If the element $x$ is such that $\|x\|= 1$, then the operator $\phi_X$ preserves the Lipschitz property, and condition~\eqref{phiProp} holds for arbitrary $m\in \RR_+$. The function $h_T(\cdot,t) = |\cdot - t|$ is extremal in inequality~\eqref{widerclass} for the real-valued functions from $H([-1,1],\RR)$. Therefore for any $x$ such that $\| x\|=1$ the function $f_t(\cdot)=|\cdot - t|\cdot x$ is extremal in the Ostrowski-type inequality for Banach space-valued functions.

In the majority results that we know (see e.g.~\cite{Anastassiou03,Anastassiou12,Kovalenko23b,Kovalenko21a}) for extremal problems on classes of non-numeric-valued functions $f\colon T\to X$,  extremal functions are built based on the real-valued extremal function for the extremal problem on the corresponding class of real-valued functions $f\colon T\to \RR$: if $f_e\colon T\to \RR$ is an extremal function in the real-valued case, then the function $f_e\cdot x$ usually becomes an extremal function in the non-numeric-valued situation for some specially chosen element $x\in X$. This corresponds to the described above approach in the case, when $\phi_X$ and $\phi_Y$ are operators of multiplication by some elements.

More generally, if for the operators $\lambda \colon H(T,M)\to M$ and $\Lambda=\lambda$, fixed $t\in T$ and $f\in H(T,M)$ inequality~\eqref{Agree} holds, 
 then  
\[
h_M (\lambda f (\cdot),\lambda f(t ))\le \lambda (h_T(\cdot,t)).
\]
If in addition, property~\eqref{phiProp} holds with $Y = M$ and $\phi_Y$ being the identity function, then the latter inequality become equality on the function  $h_T(f(\cdot,t))$. Function~\eqref{abstractOstrowskiHwExtremal} if obtained from  this function as a result of applying to it the operator $\phi_X\colon H(T,M)\to H(T,X)$.

Recall that an operator $\phi\colon H(T,M)\to X$ can be considered as an operator $\phi\colon M\to X$ (if $m\in M$, then $\phi(m):=\phi(f)$, where $f(\cdot)\equiv m$ on $T$).
\begin{corollary}\label{th::convexifyingOperator}
       Assume that operators $\Lambda, \lambda$ and $\phi_X$ and $\phi_Y=\phi_X$ satisfy the conditions of Theorem~\ref{th::generalOstrowskiHw} with $X=Y$. Let also there exist an operator $P\in H(X,X)$ such that 
    \begin{enumerate}
        \item $ \Lambda f=(\Lambda\circ P)f=(P\circ \Lambda)f \;\;\;\forall f\in H(T,X)$;
        \item $\Lambda f(t)=(\Lambda\circ P)f(t)=Pf(t)\;\;\;\forall t\in T\; \forall f\in H(T,X)$;
      % (such operator $\Lambda$ is natural to call $P$-averaging);
       \item For any $t\in T$
       \begin{equation}\label{PphiProp}
           h_X((P\circ\phi_X)(m),(P\circ\phi_X)(\theta))=m, \text { if } m=\lambda(h_T(\cdot,t)).
       \end{equation}
           \end{enumerate}
    Then for arbitrary  $t\in T$ the following sharp inequality holds:
    \begin{equation}\label{Cor}
    h_{X}(\Lambda f(\cdot), Pf(t))\le \lambda(h_T(\cdot,t)).
    \end{equation}
    The inequality becomes equality for the function
    \[    \Tilde{f}_t(\cdot)=(P\circ\phi_X)(h_X(\cdot,t)).
    \]
   If $P = {\rm Id}$ (the identity operator) satisfies the above conditions, then inequality~\eqref{Cor} has the form
    \[
    h_{X}(\Lambda f(\cdot), f(t))\le \lambda(h_T(\cdot,t)). 
    \]
    and becomes equality for function~\eqref{abstractOstrowskiHwExtremal}.
\end{corollary}

\begin{proof}
It is easy to check that operators $\Lambda, \lambda$ and $\widetilde{\phi}_X=P\circ \phi$ and $\widetilde{\phi}_Y=\widetilde{\phi}_X$ instead of $\phi_X$ and $\phi_Y$ satisfy the conditions of Theorem~\ref{th::generalOstrowskiHw} with $X=Y$.
Therefore for $h_X(\Lambda f(\cdot), \Lambda f(t))$ we obtain
\[
h_X(\Lambda f(\cdot), \Lambda f(t))\le \lambda(h_T(\cdot,t)).
\]
Due to the properties of $P$
    \[
    h_X(\Lambda f(\cdot), Pf(t))=h_X(\Lambda f(\cdot), \Lambda f(t)) \le \lambda(h_T(\cdot,t)),
    \]
and the inequality~\eqref{Cor} is proved.
    
 For the function $\Tilde{f}_t(\cdot)=(P\circ\phi_X)(h_X(\cdot,t))$ we have
     \[
h_X(\Lambda \Tilde{f}_t(\cdot), \Lambda \Tilde{f}_t(t))
=
h_X((\Lambda\circ (P\circ  \phi_X))(h_T(\cdot,t)), (\Lambda\circ (P\circ\phi_X))(h_T(t,t)))
\]
\[
\stackrel{\eqref{diagram}}{=}h_X(((P\circ \phi_X)\circ\lambda)(h_T(\cdot,t)),((P\circ\phi_X)\circ\lambda)(\theta))
\]
\[
=h_X((P\circ \phi_X)(\lambda(h_T(\cdot,t)),(P\circ\phi_X)(\lambda(\theta)))
\stackrel{\eqref{PphiProp}}{=}
\lambda(h_T(\cdot,t)).
\]
Therefore inequality~\eqref{Cor} becomes equality for the function $\Tilde{f}_t(\cdot)$.

The last statement of the Corollary is obvious.
\end{proof}

The case when $\Lambda$ is the integral operator and $P$ is the convexifying operator for multi-valued (see e.g.~\cite{nira2014}), $L$-space-valued (see e.g.~\cite{Kovalenko21a,Vahrameev}), or quasilinear-space-valued functions (see e.g.~\cite{Aseev}), is an important example of the operators that satisfy the conditions of Corollary~\ref{th::convexifyingOperator}. The case, when  $\Lambda$ is the integral operator, and $P$ is the identity operator occurs in the case of real-valued functions and functions with values in Banach spaces.
Thus in the case $M=\RR_+$ many known Ostrowski-type inequalities for real-valued, multi-valued and fuzzy-valued functions, as well as for functions with values in Banach spaces (in particular, random processes) and in $L$-spaces follow from Theorem~\ref{th::generalOstrowskiHw} and Corollary~\ref{th::convexifyingOperator} with appropriately chosen spaces  $T, X$ and  $Y$ and operators $\Lambda, \lambda, P, \phi_X$ and $\phi_Y$.

The only result that we know, where related questions were considered for $M\neq \RR_+$, is article~\cite{Kovalenko23e}.

%\newpage

\section{Classes $H^\omega(T,X)$ and inequalities of the Ostrowski type}\label{s::Hw}

\begin{definition}
Let $h_M$ be an $M$-distance in a set $M$.
A function $\omega\colon M\to M$ is called a {\it modulus of continuity}, if it satisfies the following properties:
\begin{enumerate}
    \item $\omega(\theta)=\theta$;
    \item $\omega$ is non-decreasing i.e., $\omega(m_1)\leq \omega(m_2)$, whenever $m_1\leq m_2$;
\item $\omega$ is semi-additive in the following sense: for all $m_1,m_2\in M$
$$
  h_M(\omega(m_1),\omega(m_2))
 \le 
 \omega(h_M(m_1,m_2)).
$$
\end{enumerate}
\end{definition}
In the case of  $M=\RR_+$  a modulus of continuity as an independent notion was introduced by Nikolsky~\cite{Nikolsky46}.

\begin{definition}
 Let a modulus of continuity $\omega$ and two $M$-distance spaces $(T,h_T)$, $(X,h_X)$ be given. We consider the classes
$$
    H^\omega(T,X)=\{ f\colon T\to X\colon    h_X(f(t_1),f(t_2))\le \omega(h_T(t_1,t_2))\;\;\forall  t_1, t_2\in X\}.
$$
\end{definition}

Classes  $H^\omega(T,X)$ play an important role in approximation theory. Many papers are devoted to solutions of different extremal problems for these classes. Some results for real-valued functions can be found e.g., in~\cite{ExactConstants,bagdasarov1998,Stepanets18}. Some results regarding extremal problems for classes $H^\omega(T,X)$ of functions with non-numeric values can be found in~\cite{Kovalenko23b,Kovalenko21a,Babenko15,Kovalenko22b,Drozhzhina, Kovalenko20a}.

Observe that the class $H(T,X)$ is a partial case of the class $H^\omega(T,X)$ in the case, when $\omega={\rm Id}$, where ${\rm{Id}}\colon M\to M$ is the identity mapping. On the other hand, as it is easy to see, the function  $h^\omega_T\colon T\times T\to M,$ given by the formula
\[
h^\omega_T(t_1,t_2)=\omega(h_T(t_1,t_2))
\]
is a new $M$-distance in  $T$, which becomes a $M$-pseudo metric, if  $h_T$ is an $M$-pseudo metric. Consideration of the classes  $H^\omega (T,X)$ with different $\omega$ and fixed distance $h_T$ in $T$ allows to appreciate the properties of the functions  $f\colon T\to X$ in a more detailed manner. This makes the classes   $H^\omega (T,X)$  important for approximation theory.

If in Theorem~\ref{th::generalOstrowskiHw} the $M$-pseudo metric in $T$ is understood as  $h^\omega_T$, then we obtain the following
\begin{corollary}
For arbitrary modulus of continuity  $\omega$ and arbitrary function $f\in H^\omega(T,X)$ the following inequality holds:
\[
 h_Y(\Lambda f(\cdot), \Lambda f(t))\le 
\lambda (\omega(h_T(\cdot,t))),
\]
which is sharp under the corresponding conditions and becomes equality on the function
\[
f_{\omega,t}(\cdot)=\phi_X(\omega(h_T(\cdot,t))).
\]
    
\end{corollary}

\section{On agreement of $M$-distances}\label{s::agreementOfDistances}

\begin{definition}
A partially ordered set $M$ with a smallest element $\theta$ will be called a {\it partially ordered monoid}, if an associative binary operation $+$ is defined in $M$ and the following properties hold:
\begin{enumerate}
    \item For all $m\in M$, $\theta+m=m=m+\theta.$
    \item If $m,n\in M$ are such that  $m\leq n$, then $m+p\leq n+p$ for all $p\in M$. 
\end{enumerate}
\end{definition}

\begin{definition}
An element $s$ in a partially ordered set $M$ is called a {\it supremum} of two elements $m,n\in M$, if the following two conditions are satisfied
\begin{enumerate}
    \item $s\geq m$ and $s\geq n$;
    \item If $u\geq m$ and $u\geq n$, then $u\geq s$.
\end{enumerate}
If a supremum of $m,n\in M$ exists, then we denote it by $\sup\{m,n\}$.
\end{definition}

\begin{definition}
A mapping $h_X\colon X\times X\to M$ is called an {\it $M$-metric}, if the following conditions hold:
\begin{enumerate}
    \item For all $x,y\in X$, $x=y$ if and only if $ h_{X}(x,y)=\theta$;
    \item  For all $x,y\in X$, $h_{X}(x,y)=h_{X}(y,x)$;
    \item For all $x,y,z\in X$, $h_{X}(x,y)\le h_{X}(x,z)+h_{X}(z,y)$.
\end{enumerate}
\end{definition}
%{\color{red}
In this section we give a sufficient condition on an $M$-metric $h_M$ in a partially ordered monoid $M$ to agree with an arbitrary $M$-metric  $h_X$ on a set $X$. Before doing so, we note that generally speaking $M$-metric  $h_X$ need not agree with $h_M$. For example, if $M = \RR_+$, $h_X$ is a metric such that $0 < h_X(\alpha,\beta) <1$ for some $\alpha, \beta\in X$, and $h_M$ is the discrete metric on $\RR_+$ (i.e., $h_M(a,b) = 0$, if $a = b$ and $h_M(a,b) = 1$ for all $a\neq b$), then inequality~\eqref{metricsAgreementCond} does not hold for $x = x_1 = \alpha$ and $x_2 = \beta$. Moreover, an $M$-metric $h_M$ does not necessarily agree with itself. Consider for example $M =\RR_+$, and let 
$$
h_M(a,b) = \begin{cases}
    0, & a= b = 0,\\
    \frac 34, & \text{exactly one of } a,b\text{ is } 0,\\
    \min\left\{1,\left|\ln\frac ab\right|\right\}, & a\neq 0\text{ and } b\neq 0.
\end{cases}
$$
It is easy to verify that it is actually a metric on $M$ (the fact that this function satisfies the property $h_M(r a, rb) = h_M(a,b)$ for all $a,b\geq 0$ and $r > 0$ allows to reduce the number of different cases to consider during verification of the triangle inequality). 
Since $\ln 2 < \frac 34$, for $x = x_1 = 1$, $x_2 = 2$ inequality~\eqref{metricsAgreementCond} with $h_X$ substituted by $h_M$ does not hold.
%}
 
\begin{lemma}\label{l::metricsByE}
Let $M$ be a partially ordered monoid and assume there is a function $e\colon M\times M\to M$  such that for all $x,y,z\in M$ the following properties hold:
\begin{gather*}
x\leq y\iff e(x,y) = \theta;\\
e(x,\theta) \leq  x;\\
e(x,y)\le e(x,z)+e(z,y); \notag \\
e(z+x,z+y)\le e(x,y).
\end{gather*}
If for arbitrary $x,y\in M$ the supremum $\sup\{x,y\}$ exists, then 
$$
h_M(x,y) = \sup\{e(x,y), e(y,x)\}
$$
is an $M$-metric. Moreover, arbitrary  $M$-metric $h$ agrees with $h_M$.
\end{lemma}
For example, if $M=\RR_+$, the function $e(x,y)= \max\{x-y,0\}$, $x,y\in\RR_+$, satisfies the conditions of Lemma~\ref{l::metricsByE}. In this case $h_M(x,y)=|x-y|$.

\begin{proof}
We prove that $h_M$ is an $M$-metric first. If $x\in M$, then
$$
h_M(x,x) = \sup\{e(x,x), e(x,x)\} =  \sup\{\theta, \theta\} = \theta.
$$
Moreover, if $h_M(x,y) = \theta$, then $e(x,y) = e(y,x) = \theta$, hence $x\leq y$ and $y\leq x$, thus $x = y$. Since $\sup\{a, b\} = \sup\{b, a\}$ for all $a,b\in M$, we obtain that $h_M(x,y) = h_M(y,x)$ for all $x,y\in M$.  Finally, for all $x,y,z\in X$,
\begin{gather*}
h_M(x,y) = \sup\{e(x,y),e(y,x)\}
 \le 
\sup\{e(x,z)+e(z,y), e(y,z)+e(z,x)\}
\\ \le 
\sup\{e(x,z),e(z,x)\}+\sup\{e(z,y),e(y,z)\}
=
h_M(x,z) + h_M(z,y).
\end{gather*}

Let $x\leq y$ and $z\in M$. Then $e(x,y) = \theta$, and
$$
e(y,z) = e(x,y) +  e(y,z) \geq e(x,z)
$$
i.e., the function $e$ is non-decreasing in its first variable.

Finally, if $h$ is  an $M$-metric on a set $T$, then for arbitrary $t,t_1,t_2\in T$,
\begin{gather*}
   h_M(h(t,t_1),h(t,t_2))
   =
   \sup\{ e[h(t,t_1),h(t,t_2)],e[h(t,t_2),h(t,t_1)]\}
   \\ 
   \le \sup\{ e[h(t,t_2)+h(t_2,t_1),h(t,t_2)],e[h(t,t_1)+h(t_1,t_2),h(t,t_1)]\}
   \\  
   \le \sup\{ e[h(t_2,t_1),\theta],e[h(t_1,t_2),\theta]\}
   = 
   e(h(t_1,t_2),\theta)\leq h(t_1,t_2).
\end{gather*}
\end{proof}
The idea to use an order-defining function $e$ (with properties similar to the ones stated in the lemma) as a tool to define partially ordered metric spaces was introduced in~\cite{Kovalenko20c}.

\bibliographystyle{myPlain}
\bibliography{bibliography}

\end{document}